\documentclass[11pt]{amsproc}

\newtheorem{Theorem}{Theorem}

\theoremstyle{definition}

\theoremstyle{remark}

\numberwithin{equation}{section}

\newcommand{\abs}[1]{\lvert#1\rvert}

\newcommand{\Z}{{\mathbb Z}}

\begin{document}

\title{Two remarks on $C^\infty$ Anosov diffeomorphisms}


\author{Shigenori Matsumoto}
\address{Department of Mathematics, College of
Science and Technology, Nihon University, 1-8-14 Kanda, Surugadai,
Chiyoda-ku, Tokyo, 101-8308 Japan
}
\email{matsumo@math.cst.nihon-u.ac.jp
}
\thanks{The author is partially supported by Grant-in-Aid for
Scientific Research (C) No.\ 20540096.}


\date{\today }

\begin{abstract}
Let $M$ be a closed oriented $C^\infty$ manifold and 
$f$ a $C^\infty$ Anosov diffeomorphism on $M$.
We show that if $M$ is the two torus $T^2$, then
$f$ is conjugate to a hyperbolic automorphism of $T^2$,
either by a $C^\infty$ diffeomorphism or by a singular homeomorphism.
We also show that for general $M$, if $f$ admits
 an absolutely continuous invariant measure $\mu$,
 then $\mu$ is a $C^\infty$ volume.
The proofs are concatenations of well known results in the field.
\end{abstract}

\maketitle

\section{Conjugacy}

Let $f$ be a $C^\infty$ Anosov diffeomorphism on the two torus $T^2$.
Then
$$
A=f_*\in {\rm Aut}(H_1(T^2,\Z))={\rm SL}(2,\Z)
$$
defines a hyperbolic automorphism of the abelian Lie group $T^2$, and
$f$ is isotopic to $A$.
It is known \cite{F,M} that $f$ is conjugate to $A$
by a homeomorphism $h$ which
is isotopic to the identity: $h\circ A=f\circ h$. It is well known that the conjugacy
$h$ is a bi-H\"older homeomorphism. Also it is easy to show that
$h$ is unique. 
Let us denote by $m$ the normalized Haar measure of $T^2$. A homeomorphism $h$
of $T^2$ is said to be {\em singular} if there is an $m$-conull Borel
set $E$ such that $h(E)$ is $m$-null.
Our first result is the following.

\begin{Theorem} \label{t1}
The conjugacy $h$ is either a $C^\infty$ diffeomorphism or a
singular homeomorphism.
\end{Theorem}

\medskip

{\sc Proof}.
Let  $TT^2=E^u\oplus E^s$ be the hyperbolic splitting
associated with $f$.
By a dimensional reason, it is a $C^1$ splitting \cite{MM1}.
Fix a translation invariant $C^\infty$ Riemannian metric $g$ on $T^2$.
The derivative of $f$ along $E^u$ (resp.\ $E^s$)
measured with respect to $g$ is denoted by $J^sf$
(resp.\ $J^uf$). These are $C^1$ functions. 
The Gibbs measure \cite{S,B} for the potential $-\log \abs{J^uf}$
(resp.\ $\log \abs{J^sf})$ is denoted by $\mu_+$ (resp.\ $\mu_-$).

Let $f$, $A$ and $h$ be as above. First consider
the case where $f$ does not admit an a.\ c.\ i.\ m.\
(absolutely continuous invariant measure). Then the
$f$-invariant measure $h_*m$ 
is singular to $m$.

To show this, notice that $h_*m$ is decomposed into two parts; one absolutely continuous
and the other singular. Since $h$ is a $C^\infty$ diffeomorphism, it
leaves each part invariant.  But the absolutely continuous part must
be zero since by the assumption there is no a.\ c.\ i.\ m.\ for $f$.

Thus $h$ maps the measure $m$ to a singular measure $h_*m$.
Since $h_*m$ is singular,
there is an $m$-null set $E'$ such that $h_*m(E')=m(h^{-1}E')=1$.
The conjugacy $h$ maps the $m$-conull set $h^{-1}(E')$ to the $m$-null set
$E'$, showing that  $h$ is a singular homeomorphism.

Next consider the case where $f$ admits an a.\ c.\ i.\ m.\ $\mu$.
Then we have $\mu=\mu_+=\mu_-$ (Proof of Corollary 1 of
\cite{S}, Corollary 4.13 of \cite{B}).
In particular an a.\ c.\ i.\ m.\ $\mu$ is unique and
ergodic. The induced measure $h_*m$ is also
ergodic. Therefore either $\mu$ and $h_*m$ are mutually singular
or coincide. In the former case, we argue just as before,
to conclude that the conjugacy $h$ is a singular homeomorphism.

Finally assume that $\mu_+=\mu_-=h_*m$. These are the Gibbs measures
of three potentials, $-\log \abs{J^uf}$, $\log \abs{J^sf}$ and a constant.
By Section 3.4 of \cite{S}, these three functions, with the identical
Gibbs measure, are
mutually cohomologous modulo constant. That is, there are
continuous functions $v_1$, $v_2$ and constants $c_1$, $c_2$
such that
$$-\log \abs{J^uf}=v_1\circ f-v_1+c_1,
$$
$$\log \abs{J^sf}=v_2\circ f-v_2+c_2.$$

This shows that the Lyapunov exponents of all periodic orbits are
the same. By Theorem 1 of \cite{MM2}, the conjugacy $h$ is
a $C^\infty$ diffeomorphism.
The proof of Theorem 1 is complete.
\qed

\section{Absolutely continuous invariant measure}

Let $M$ be a closed oriented $n$-dimensional $C^\infty$ manifold
and $f$ a $C^\infty$ Anosov diffeomorphism on $M$.
Let $g$ be a $C^\infty$ Riemannian metric on $M$, 
and $m$ the normalized measure given by the volume
form associated with  $g$.

\begin{Theorem}
Assume $f$ admits an a.\ c.\ i.\ m.\ $\mu$ with
density $\varphi$: $\mu=\varphi m$, $\varphi\in L^1(m)$.
Then the density $\varphi$ is a positive $C^\infty$ function.
\end{Theorem}

{\sc Proof}.
Let  $TM=E^u\oplus E^s$ be the hyperbolic splitting
associated with $f$.
Denote the Jacobian along $E^u$ (resp.\ $E^s$)
measured with respect to $g$ by $J^uf$ (resp.\ $J^sf$).
The total Jacobian measured with respect to $g$
is denoted by $Jf$.
All these are continuous real valued functions on $M$.

Define another continuous Riemannian metric 
$g'$ by $g'=g\vert_{E^u}\oplus g\vert_{E^s}$.
Thus $E^u$ and $E^s$ are perpendicular with respect to $g'$.
Let $m'$ be the normalized measure given by the volume
form associated with $g'$.
We have $m'=e^am$ for a continuous function $a$.

Denote by $J'f$ the total Jacobian with respect to $g'$.
Then we have
\begin{equation} \label{e11}
\log\abs{J'f}=\log\abs{J^uf}+\log\abs{J^sf}.
\end{equation}

By \cite{S,B}, we have $\mu=\mu_+=\mu_-$,
where $\mu_+$ (resp.\ $\mu_-$) is the Gibbs measure for
the potential $-\log \abs{J^uf}$ (resp.\ $\log\abs{J^sf}$).
Then by \cite{S}, 
$\log \abs{J^uf}+\log\abs{J^sf}$ is cohomologous to a constant.
Thus by (\ref{e11}), we have
\begin{equation} \label{e12}
\log\abs{J'f}=b\circ f- b+C.
\end{equation}
for a continuous function $b$ and a constant $C$.

On the other hand, by the invariance of the a.\ c.\ i.\ m.\ 
$\mu=\varphi e^{-a}m'$,
we have $\mu$-almost everywhere
\begin{equation} \label{e13}
\log\abs{J'f}=(a-\log\varphi)\circ f-(a-\log\varphi)
\end{equation}

Now by (\ref{e13}), we have $\mu(\log\abs{J'f})=0$.
This implies that $C=0$ in (\ref{e12}).
Then (\ref{e12}) implies the invariance of the
measure $e^{-b}m'=e^{-b+a}m$. Moreover, adding an
appropriate constant to $b$,
we may assume that $e^{-b+a}m$ is a probability measure.
By the uniqueness of
the a.\ c.\ i.\ m., we have $\mu=e^{-b+a}m$.
That is, the density of $\mu$ is positive and continuous.
Then by Corollary 2.1 of
\cite{LMM}, we obtain that  $e^{-b+a}$ is a $C^\infty$
function. The proof is complete.
\qed

\end{document}